\documentclass[10pt]{amsart}

\title[Dirichlet Problem for Curve Shortening]{The Dirichlet problem for curve shortening flow} 
\date{\today}
\author{
Paul T.~Allen
}
\address{
Department of Mathematical Sciences,
Lewis \& Clark College,
Portland, U.S.A.
}

\author{
Adam Layne
}
\address{
Department of Mathematics,
University of Oregon,
Eugene, U.S.A.
}

\author{
Katharine Tsukahara
}
\address{
Indiana University School of Medicine,
Indianapolis, U.S.A.
}

 	\usepackage{setspace}
	\usepackage{indentfirst} 
	\usepackage{amsmath}
	\usepackage{amssymb}
	\usepackage[usenames,dvipsnames]{color}
	\usepackage{amsfonts}
	\usepackage{verbatim}
	\usepackage[mathscr]{eucal}
	\usepackage{amsthm}
	\usepackage[numbers]{natbib} 

	\usepackage{ mathrsfs }
	
	\usepackage{hyperref}

\theoremstyle{plain} 
\newtheorem{lemma}{Lemma}[section]
\newtheorem{prop}{Proposition}[section]
\newtheorem{theorem}{Theorem}

\newcommand{\scL}{\mathcal{L}}
\newcommand{\scrL}{\mathcal{L}}
\newcommand{\bbR}{\mathbb{R}}
\newcommand{\bbH}{\mathbb{H}}
\newcommand{\bbS}{\mathbb{S}}
\newcommand{\scrS}{\mathscr{S}}

\begin{document}

\maketitle

\begin{abstract}

We investigate the evolution of curves with fixed endpoints under the curve shortening flow, which evolves curves in proportion to their curvature.
Curve shortening flow has a long history, and closed curves and curves in the Euclidean plane have been extensively studied.
Using a distance comparison by Huisken, we determine the long-term behavior of open curves with fixed endpoints evolving in certain convex domains on surfaces of constant curvature.
Specifically, we show that such curves do not develop singularities, and evolve to geodesics.

\end{abstract}

\section{Introduction}
We consider the Dirichlet problem for the curve shortening flow on Riemannian surface $(\Sigma, g)$.
For curve $I=[a,b] \to \Sigma$ given by $u\mapsto F(u)$, the induced arc-length is given by $ds^2 = \langle \partial_u F, \partial_u F\rangle_g du^2$.
The first variation of length $\scrL= \int_I ds$ is given by
\begin{equation}
\delta\scrL[F](\xi) = - \int_I \langle K, \xi\rangle_g ds. 
\end{equation}
Here the curvature vector given by $K = \nabla_TT$, with $T = \frac{d}{ds}F$ being  the unit tangent vector.
The curve shortening flow for function $F:I\times[0,t^*)\to \Sigma$ is the gradient flow associated to $\scrL$, and is given by equating the velocity vector $V=\partial_tF$ and the curvature vector $K$:
\begin{equation}\label{main}
\partial_t F = K.
\end{equation}
We consider the initial boundary value problem associated to \eqref{main}, requiring for all $t$ that 
\begin{equation}\label{DirichletBC}
F(a,t)=A\text{ and }F(b,t)=B
\end{equation}
for fixed $A,B\in\Sigma$, and under natural hypotheses show that the flow converges to a geodesic joining $A$ and $B$.

The evolution of closed curves under \eqref{main} has received considerable attention: for planar curves see \cite{GageHamilton86},\cite{Grayson87} and also \cite{Angenent91},\cite{Hamilton95},\cite{Huisken98} among others; curves in Euclidean $3$-space are considered in \cite{Altschuler91}.
As is typical of parabolic geometric flows, singularity formation is associated with self-similar behavior; such solutions are studied in \cite{AbreschLanger86},\cite{Halldorsson2010}.
Closed curves on Riemannian surfaces are studied in \cite{Gage90b},\cite{Grayson89}; a broad class of evolutions of closed curves on surfaces is studied in \cite{Angenent90I,Angenent91II}.

To our knowledge, the only results for curves with fixed endpoints appear in \cite{Huisken98} where it is shown that embedded curves bounded by two parallel lines, and with one endpoint on each line, flow towards the unique geodesic joining the endpoints.
Our work here revisits this result and extends it to curves on surfaces, which for simplicity we assume have constant curvature and, as one may lift the problem to the universal cover, trivial topology.

Before stating our main theorem, let us carefully describe the rationale for the geometric situation we consider.
The compatibility condition associated to the Dirichlet boundary condition (see Section \ref{S-Continuation} below) is that the curvature, and its derivatives, vanish at the endpoints; this is thus required of the initial curve.

We furthermore assume that the endpoints $A,B$ of the initial curve  lie on `barrier' geodesics $\Gamma_a ,\Gamma_b$, and that interior of the initial curve lies in an open, convex region, which is necessarily `between' the geodesics.
On the sphere $\bbS^2$ this implies that the initial curve lies on a closed hemisphere with at most only the endpoints lying on the equator.
The reason for this condition is two-fold.
First, a maximum principle of Angenent (Proposition \ref{P:AngMax}) implies that geodesics are barriers for solutions to curve shortening flow; thus this condition ensures that the flow remains in a bounded region containing the desired limiting geodesic.
(As one might expect, uniqueness issues arise when $A,B$ are conjugate; see Section \ref{S-Converge}.)

Secondly,  if self-intersections form, one expects loops to contract into a cusp singularity; see \cite{Angenent91}.
The same maximum principle implies that interior self-intersections cannot form; for closed curves this ensures that initially embedded curves remain so.
In our situation, however, the absence of a barrier geodesic permits initial curves which form self-intersections by means of the curve passing over the endpoint.
The geodesics $\Gamma_a,\Gamma_b$ thus ensure that curves which are initially embedded remain so.

We now state our main result.
Let $\Sigma$ be be one of $\bbS^2,\bbR^2,\bbH^2$ with Gaussian curvature $S = 1,0,-1$, respectively. 
\begin{theorem}\label{MainTheorem}
Let $F_0:[a,b]\to\Sigma$ be an embedded curve lying in a closed, convex region ${\Omega}$ with interior lying in $\text{\emph{int}}(\Omega)$.
Furthermore, suppose that the boundary of $\Omega$ consists of two (possibly coincident) geodesics $\Gamma_a,\Gamma_b$ and that $A=F_0(a)\in \Gamma_a$, $B=F_0(b)\in\Gamma_b$ but $A\ne B$.
Finally, suppose that the curvature of the curve given by $F_0$  and all its derivatives vanish at $A,B$.

Then there exists $F:I\times[0,\infty)\to\Sigma$, with $F(\cdot,0) = F_0$, satisfying \eqref{main}-\eqref{DirichletBC} such that the curvature of $F(I,t)$ tends smoothly to zero as $t\to\infty$.
If $A,B$ are not conjugate, then $F(I,t)$ tends to the unique geodesic joining the endpoints.
\end{theorem}

In the case of $\Sigma = \bbR^2$, the result is due to Huisken \cite{Huisken98}, and we follow his approach here.
In Section \ref{S-Preliminary} we record a number of primary statements concerning curve shortening flow on surfaces.
Section \ref{S-EndpointControl} makes use of the aforementioned maximum principle to establish control of the flow near the endpoints, ensuring that the flow remains embedded.
Finiteness of the curvature is sufficient to ensure that the flow may be continued; we record this in Section \ref{S-Continuation}.
In lieu of directly controlling the curvature, we follow the approach in \cite{Huisken98} and establish, in Section \ref{S-DistanceComparison}, a lower bound on the ratio of extrinsic and intrinsic distances.
The estimate has the advantage of being scale-invariant.
We record in Section \ref{S-Monotonicity} the monotonicity formula of \cite{Huisken90}.
The argument now proceeds by contradiction; we assume that a finite-time singularity ceases the flow.
In Section \ref{S-Singularities} we prove general properties of singularities before using the distance comparison estimate to rule out the possibility that they occur; this is accomplished by means of blowup analyses in Sections \ref{S-TypeOneSingularities}, \ref{S-TypeTwoSingularities}.
Finally, we discuss long-time convergence in Section \ref{S-Converge}.

We remark that many of the arguments here can be extended to the case of curves in a general Riemannian surface $(\Sigma,g)$ with pointwise bounded curvature.
We have primarily used the assumption that $\Sigma$ has constant curvature to simply computations, most notably in Lemma \ref{SecondVariationD} and in our application of \eqref{Wirtinger} in Section \ref{S-Converge}.
We have also chosen to work in the smooth setting, which necessitates the higher-order compatibility condition that all derivatives of the curvature vanish at the endpoints; see section \ref{S-Continuation}.
 
We thank J.~S.~Ely for providing numerical simulations, I.~Stavrov for helpful conversations, F.~Morgan for a helpful comment.
Partial funding for this work comes from the J.S.~Rogers Science Research Program at Lewis \& Clark College.

\section{Preliminary computations}\label{S-Preliminary}

We record here several basic facts concerning the geometry of solutions to \eqref{main}, which are recorded in \cite{GageHamilton86} for $\Sigma = \bbR^2$ and in \cite{Gage90b},\cite{Grayson89} for general surfaces.

Setting $\kappa = \langle K,N\rangle$, for some choice of unit normal vector $N$ along $F(I,t)$, the Frenet equations are
\begin{equation}\label{frenet}
\nabla_T T  = \kappa N \quad \text{and} \quad \nabla_T N =- \kappa T.
\end{equation}
Here $\nabla$ is the connection associated to $g$.

Let $v = \langle \partial_uF,\partial_uF\rangle^{1/2}$; direct computation shows that if \eqref{main} is satisfied, then
\begin{equation}	\label{E:v-evolution}
\partial_tv = - \kappa^2 v,
\end{equation}
which implies the following two lemmas.

\begin{lemma}	\label{L-evol}
The evolution formula for the length of the curve under curve shortening flow is given by
\begin{equation}
\frac{d}{dt}\scL = -\int_I \kappa^2 ds.
\end{equation}
\end{lemma}
Note that $\scrL$ is bounded below by the geodesic distance between the fixed endpoints.

\begin{lemma}\label{ds-commute}
The commutator of  $V=\partial_tF$ and $T=\partial_sF$ under curve shortening flow is given by
\begin{equation}
\left[ V,T \right] = \kappa^2 T.
\end{equation}
\end{lemma}

Lemma \ref{ds-commute}, together with \eqref{frenet}, implies the following:

\begin{lemma}\label{TN-evol}
The evolutions of $T$ and $N$ under curve shortening flow are given by
\begin{equation}
\nabla_{V} T = (\partial_s\kappa) N 
\quad \text{and} \quad 
\nabla_{V} N = -(\partial_s\kappa) T.
\end{equation}
\end{lemma}
Since $V = \partial_tF = K = \kappa N$ along the flow, the definition of the curvature tensor $R$ associated to $g$ yields
\begin{equation}
\nabla_{V}\nabla_{T} - \nabla_T\nabla_{V} = \kappa^2\nabla_T + \kappa R(T,N).
\end{equation}
This in turn implies the following formula for the evolution of $\kappa$.
\begin{lemma}\label{L-EvolutionOfCurvature}
Under the curve shortening flow we have
\begin{equation}\label{EvolveKappa}
\partial_t\kappa = \partial_s^2\kappa + \kappa^3 + S\kappa,
\end{equation}
where $S$ is the Gaussian curvature of $\Sigma$ (precomposed with $F$).
\end{lemma}

The parabolic equation \eqref{EvolveKappa} is central to the study of solutions to the curve shortening flow, and is used to establish geometric control in a variety of ways.
We record the following from \cite{Altschuler91} (see also \cite{Grayson87}), indicating that the total turning angle is non-increasing.
\begin{lemma}\label{TurningAngleEstimate}
Under the curve shortening flow we have
\begin{equation}
\frac{d}{dt}\int_I |\kappa|ds = -2\sum_{\{u \mid \kappa(u)=0\}}|\partial_s\kappa(u)|.
\end{equation}
\end{lemma}

\section{Control at endpoints}\label{S-EndpointControl}

Results of Angenent \cite{Angenent91, Angenent91b} provide further geometric control over solutions to the curve shortening flow.
We recall from \cite{Angenent91} the following.
\begin{prop}\label{P:AngMax} 
Suppose $u(x,t)$ is continuous on $[x_0,x_1]\times[t_0,t_1]$ and is a classical solution of 
\begin{equation}\label{Angenent-main}
\partial_tu = a\partial_x^2u+b\partial_xu+cu
\end{equation}
on $(x_0,x_1)\times(0,T]$.
Here $a>0$ is $C^2$, $b$ is $C^1$, and $c$ is $C^0$ on $[x_0,x_1]\times[t_0,t_1]$.
We assume $u(x_i,t)\neq 0$ for $i=1,2$ and $t\in [t_0,t_1]$.

Let $z(t)$ be the number of elements in the zero set $Z(t)=\{x \in [x_0,x_1]  \mid u(t,x) = 0\}$.
For $t\in (0,T]$, we have that $z(t)$ is finite.

Furthermore, if for some $t_c \in(0,T)$, we have an degenerate inflection point (i.e.,~$u(x_c,t_c) = \partial_xu(x_c,t_c) = 0)$, then for all $t_- < t_c \leq t_+$ we have $z(t_-) > z(t_+)$.
\end{prop}

The proposition has several important consequences. Locally, solutions to curve shortening flow can be written as a graph, in which case \eqref{main} takes the form \eqref{Angenent-main}. 
Proposition \ref{P:AngMax} yields the following.
\begin{lemma}
The number of interior intersections of solutions to the curve shortening flow cannot increase.
\end{lemma}
In particular this implies that no interior self-intersections form.
It also implies that geodesics, which are stationary solutions of curve shortening flow, serve as barriers for solutions to the curve shortening flow. The existence of these barriers implies the following.
\begin{lemma}\label{L-RemainsBounded}
Under our hypotheses, solutions to the curve shortening flow remain in a bounded, convex region of $\Sigma$.
\end{lemma}

By considering a foliation of a neighborhood of the endpoints $A$, $B$ by geodesics, we obtain a gradient estimate near the endpoints.
\begin{lemma}
\label{L-GradientEstimate}
There exists an $\epsilon$-neighborhood of $\Gamma_a$ in which the curves $F(I,t)$, when viewed as graphs over a fixed geodesic $\Gamma_a^*$ passing through $A$ and transverse to $\Gamma_a$, have uniformly bounded gradient. 
\end{lemma}
The same statement holds in a neighborhood of $\Gamma_b$. 
The gradient bound implies that the curvature near the endpoints is uniformly bounded.
\begin{lemma}
\label{L-BoundaryCurvatureBound}
There exists an $\epsilon$-neighborhood of $\Gamma_a$ and of $\Gamma_b$ in which $\kappa$ is uniformly bounded.
\end{lemma}

Finally, we have the following.
\begin{lemma}
Under our hypotheses, the solution to the curve shortening flow remains embedded so long as the flow exists.
\end{lemma}

\section{Local existence \& continuation}\label{S-Continuation}

The system \eqref{main} is not strictly parabolic; there are zeros present in the symbol, associated to diffeomorphism invariance of the problem.
Short-time existence for the the curve shortening flow is discussed in \cite{GageHamilton86}, where the Nash-Moser inverse function theorem  is used; see \cite{Hamilton82b},\cite{Angenent90I}. 

In the current case, it suffices to choose coordinates so that the curve is a graph of height $f$ over the initial curve.
In these coordinates, the evolution of $f$ is given by
\begin{align}
f_t = \kappa \frac{1}{\langle \partial_0,N \rangle_g}
\end{align}
where $x^0$ is the height coordinate.
This is (for small $t$) a parabolic equation on $I$ for $f$ which vanishes at the parabolic boundary.
Existence and uniqueness theorems for equations of this type follow from classical theory; see Theorems 4.1 and 1.1 in Chapter VI of \cite{MR0241822}, which imply existence of a smooth solution.
The smoothness of the solution is contingent upon higher-order compatibility conditions being satisfied at the endpoints.
These conditions are implied by the vanishing of the the curvature and its derivatives at the boundary at the initial time.
The problem remains strictly parabolic as long as the curvature of the solution is bounded, and so we have the following existence theorem.

\begin{prop}\label{ContinuationCriterion}
If the curvature $S$ of $(\Sigma,g)$ is bounded, and the curvature and its spatial derivatives vanish on the boundary of the interval, then there exists for short time a solution $F(t):I\to\Sigma$ of the curve shortening flow \eqref{main}.
Furthermore the solution continues to exist, and is smooth, so long as $\kappa(t)$ remains bounded.
\end{prop}

Following \cite{Grayson87},\cite{Grayson89}, Lemmas \ref{ds-commute} and \ref{L-EvolutionOfCurvature} can be used to establish (schematic) formulas for the evolution of derivatives of $\kappa$.
\begin{lemma}\label{L-EvolutionOfDerivativesOfCurvature}
The derivatives of curvature evolve according to
\begin{equation}\label{HowDerivativesEvolve}
\partial_t\left[ \partial^n_s\kappa\right] 
= \partial_s^2\left[\partial_s^n\kappa\right] 
+ \left((n+3)\kappa^2 + S\right) \partial_s^n\kappa
+ P_n + F_n,
\end{equation}
where $P_n$ is a polynomial in $\kappa, \partial_s\kappa, \dots, \partial_s^{n-1}\kappa$ for which each monomial contains at most $n-1$ derivatives and $F_n = \sum_{i=1}^n {n \choose i} (\partial_s^{n-i} \kappa)( \partial_s^i S)$.
\end{lemma}

Inductively making use of the maximum principle for parabolic equations we conclude that bounds on $\kappa$ imply bounds on derivatives of $\kappa$.
\begin{lemma}
\label{L-KappaBoundIsEnough}
Suppose  $\kappa$ is bounded on $[0,t)$. Then all derivatives $\partial_s^k\kappa$ are also bounded on $[0,t)$.
\end{lemma}

\section{Distance comparison}\label{S-DistanceComparison}
In order to establish the main theorem, we argue that no finite-time singularity occurs.
We proceed by contradiction, assuming that there is some maximal time of existence $t^*<\infty$ at which a singularity occurs.
The desired contradiction is obtained via a distance comparison estimate, established in the planar case by Huisken \cite{Huisken98}; the estimate serves as a scale-invariant proxy for a curvature bound.

For simplicity we restrict attention to the case where $(\Sigma,g)$ is one of the round sphere $\bbS^2$, the Euclidean plane $\bbR^2$, or the hyperbolic plane $\bbH^2$, and make use of geodesic polar coordinates $(\phi,\theta)$ in which the metric $g$ takes the form $g = d\phi^2 + \scrS(\phi)^2\,d\theta^2$.
Here $\scrS(\phi) = \sin\phi$ when $\Sigma = \bbS^2$, $\scrS(\phi) = \phi$ when $\Sigma = \bbR^2$, and $\scrS(\phi)=\sinh\phi$ when $\Sigma = \bbH^2$.
Recall that Lemma \ref{L-RemainsBounded} implies that solutions to the curve shortening flow remain in a bounded, convex region of $\Sigma$. 
This region is covered by such polar coordinates, centered at any point in the region.

Let $D(P,Q)$ be the geodesic distance between $P,Q\in\Sigma$.
This gives rise to the function $D(p,q;t) = D(F(p,t),F(q,t))$, defined on $I\times I\times [0,t^*)$, which measures the extrinsic distance between points on the curve given by $F(\cdot, t)$.
This can be compared to the intrinsic distance, given by 
\begin{equation}
L(p,q;t) = \int_p^q v(u,t)\,du.
\end{equation}
Without loss of generality we restrict attention to $p\leq q$.

Clearly $(D/L)(p,q;t) \leq 1$ with equality when $F([p,q], t)$ traces out a geodesic.
Furthermore, $(D/L)(p, q; t) =1$ when $p=q$.
We establish a lower bound on $D/L$ using the following interior estimate.
\begin{prop}\label{Prop-InteriorEstimate}
Let $F:I\times[0,t^*)\to\Sigma$ be a smooth, embedded solution to \eqref{main}.
Suppose at some time $\hat t\in (0,t^*)$ that $(D/L)(\cdot,\cdot\,;\hat t\,)$ attains a local minimum at $(p,q)$ in the interior of $I\times I$.
Then
\begin{equation}
\frac{d}{dt}\left[ \frac D L(p, q;t)\right]_{t=\hat t} \geq 0.
\end{equation}
\end{prop}

Supposing that $(D/L)(\cdot,\cdot\,;\hat t\,)$ attains a local minimum at $(p,q)$, let $W$ be the unit vector tangent to the geodesic from $F(p,\hat t)$ to $F(q,\hat t)$, let $V$ be a unit vector orthogonal to $W$.
For any vector $\xi$ tangent to $F(I,\hat t)$ we use the notation $\xi_p = \xi \big|_{F(p,\hat t)}$.

We now compute the variations of $D$ and $L$ along the curve $F(\cdot, \hat t)$.
The following are immediate.
\begin{lemma}\label{FirstVariationOfL}
The first variation of $L$ is given by
\begin{equation}
\delta L(p,q;\hat t)[\xi_p,\xi_q] = \langle T_q, \xi_q\rangle_g - \langle T_p,\xi_p\rangle_g,
\end{equation}
while the second variation $\delta^2 L(p,q;\hat t)(\xi_p,\xi_q) \equiv 0$.
\end{lemma}
\begin{lemma}\label{FirstVariationOfD}
The first variation of $D$ is given by
\begin{equation}
\delta D(p,q;\hat t)[\xi_p,\xi_q] = \langle W_q, \xi_q\rangle_g - \langle W_p,\xi_p\rangle_g.
\end{equation}
\end{lemma}

The first variation of $D/L$ may be computed using the formulas above.
Since $(p, q)$ is a local minimum, we have
\begin{equation}
0
= \frac{1}{\hat L}\left(\langle W_q, \xi_q\rangle_g - \langle W_p,\xi_p\rangle_g\right) 
- \frac{\hat D}{\hat L^2} \left(\langle T_q, \xi_q\rangle_g - \langle T_p,\xi_p\rangle_g\right),
\end{equation}
where we have set $\hat D = D(p,q;\hat t)$ and $\hat L = L(p,q;\hat t)$.
Taking $\xi$ such that $\xi_p=T_p$ and $\xi_q=0$, and vice versa, we obtain the following.
\begin{lemma}\label{TdotW}
\begin{equation}
\langle  W_{p},T_{p} \rangle_g 
=  \langle  W_{q},T_{q} \rangle_g 
= \frac{\hat D}{\hat L}.
\end{equation}
\end{lemma}
We now choose $\xi$ such that $\xi_p = T_p$ and $\xi_q = \epsilon T_q$, where $\epsilon = \pm1$. 
By Lemmas \ref{FirstVariationOfL}, \ref{FirstVariationOfD}, \ref{TdotW}  we have
\begin{equation}
\begin{aligned}
\widehat{\delta D} 
&= \delta D(p,q;\hat t)[ T_p,\epsilon T_q] 
= (\epsilon-1)\frac{\hat D}{\hat L},
\\
\widehat{\delta L} 
&=\delta L(p,q;\hat t)[ T_p,\epsilon T_q] 
= (\epsilon-1),
\end{aligned}
\end{equation}
and thus
\begin{equation}
\begin{aligned}
\delta^2 \left(\frac D L \right)(p,q;\hat t)[ T_p,\epsilon T_q] 
&=\frac{1}{\hat L}\widehat{\delta^2D},
\end{aligned}
\end{equation}
where $\widehat{\delta^2D} = \delta^2 D(p,q;\hat t)[ T_p,\epsilon T_q]$.
A somewhat lengthy computation using the law of cosines
yields the following formula for $\widehat{\delta^2D}$.

\begin{lemma}\label{SecondVariationD}
The second variation of $D$, with respect to $\xi$ as described above, is given by
\begin{equation}
\widehat{\delta^2D} = \langle K_q,W_q\rangle_g - \langle K_p,W_p\rangle_g  
-\frac{\scrS(\hat D/2)}{2 \scrS^\prime(\hat D/2)}\left(\epsilon\langle T_q,V_q\rangle_g - \langle T_p,V_p\rangle_g \right)^2.
\end{equation}
\end{lemma}

Lemma \ref{TdotW} implies that we may choose $\epsilon$ such that $\langle T_q,V_q\rangle_g = \epsilon \langle T_p,V_p\rangle_g$.
Applying this to the minimizing condition $\delta^2 (D/L)\geq 0$ leads to the following.
\begin{lemma}
If $(p,q;\hat t)$ is an interior local minimum for $D/L$, then $\langle K_q,W_q\rangle_g - \langle K_p,W_p\rangle_g\geq 0$.
\end{lemma}
Proposition \ref{Prop-InteriorEstimate} follows directly from this last lemma and Lemma \ref{L-evol}.
In fact, we have
\begin{equation}
\frac{d}{dt}\left[ \frac D L(p, q;t)\right]_{t=\hat t} \geq \int_{s(p)}^{s(q)} \kappa(s,\hat t)^2 ds.
\end{equation}

We now argue that the ratio $D/L$ is bounded below near the boundary of $I\times I$.
Note that by Lemma \ref{L-evol} we have $L$ bounded above which implies $D/L$ is bounded below away from the diagonal of $I\times I$.
When both $p$ and $q$ are near the same endpoint of $I$ the gradient estimate of Lemma \ref{L-GradientEstimate} implies that $L$ is bounded above by a constant multiple of $D$.
The ensuing lower bound on $D/L$ near the boundary of $I\times I$, together with the previous proposition, implies the following.

\begin{prop}\label{DistanceComparisonBound}
Let $\Sigma$ be one of $S^2$, $\bbR^2$, $\bbH^2$ and let $F$ be a solution to \eqref{main}-\eqref{DirichletBC}.
Then there exists $\theta>0$ such that $D/L\geq\theta$ on $I\times I\times [0,t^*)$.
\end{prop}

\section{Monotonicity formula}\label{S-Monotonicity}
Following \cite{Huisken90}, we establish a monotonicity formula in the case that $\Sigma$ has constant curvature.
Fix $P^*$ within the convex region of $\Sigma$ containing the curve and set $\rho(u,t) = D(P^*, F(u,t))$.
We make use of the backwards heat kernel $\Lambda:I\times[0,t^*)\to \bbR$, given by
\begin{equation}
\Lambda(u,t)= \frac{1}{ \sqrt{4 \pi ({t^*}-t)} } \exp{\left( -\frac{ \rho(u,t)^2}{4 (t^*-t)} \right)},
\end{equation}
to define
\begin{equation}
Q(t)= \int_I \Lambda(u,t)\, ds(u).
\end{equation}
Direct computation yields 
\begin{equation}\label{DerivativesOfLambda}
\begin{aligned}
\partial_t\Lambda &=	- \Lambda \left( \left( {\frac{\rho}{2(t^*-t)}} \right)^2- \frac {1}{2(t^*-t)} 
	+ \left\langle K, \frac{\rho}{2(t^*-t)}W\right\rangle_g \right) ,
\\
\partial_s\Lambda &= -\Lambda    \left\langle T ,{\frac{\rho}{2(t^*-t)}} W \right\rangle_g  ,
\end{aligned}
\end{equation}
where $W=W(u,t)$ is the unit tangent vector of the geodesic from $P^*$ to $F(u,t)$, evaluated at $F(u,t)$.
Further computation and integration by parts leads to the following formula.
\begin{lemma}
If $F$ is a solution to the curve shortening flow, then $Q$ evolves according to
\begin{multline}
\frac{d}{dt} Q(t)	
=\int_I \frac{\Lambda}{2(t^*-t)} ds - \int_I \Lambda \left | K+{\frac{\rho}{2(t^*-t)}} W  \right |^2_g  ds
\\
+\left. \left\langle \Lambda {\frac{\rho}{2(t^*-t)}} W  , T \right\rangle_g \right|_{\partial I}
 - \int_I \left\langle  \nabla_T \left(\Lambda {\frac{\rho}{2(t^*-t)}} W \right) , T  \right\rangle_g ds.
\end{multline}
\end{lemma}
We now compute $\left\langle  \nabla_T \left(\Lambda {\frac{\rho}{2(t^*-t)}} W \right) , T  \right\rangle_g$ at some point $Q^*$, making use of geodesic polar coordinates centered at $P^*$.
Without loss of generality we may assume $W^\phi=1$ and $W^\theta=0$ at $Q^*$, in which case $\langle T,W\rangle_g = T^\phi$ and thus
\begin{equation}
\left\langle\nabla_T W,T\right\rangle_g = \frac{1}{\scrS} \frac{d\scrS}{d\phi} \left(1- \left\langle T,W \right\rangle_g^2\right).
\end{equation}
Hence
\begin{multline}
\left\langle \nabla_T	\left( \Lambda {\frac{\rho}{2(t^*-t)}} W \right),T \right\rangle_g
=\left( \frac{1}{2(t^*-t)}- \left({\frac{\rho}{2(t^*-t)}}\right)^2 \left\langle T ,W \right\rangle_g^2\right) \Lambda
\\
+ \frac{\Lambda}{2(t^*-t)} \left(1- \left\langle T,W \right\rangle_g^2\right) \left(\rho\frac{\scrS'(\rho)}{\scrS(\rho)}-1\right) .
\end{multline}

From this we obtain a version of Huisken's monotonicity formula (compare \cite{Huisken90}).
\begin{lemma}
Under the curve shortening flow, the quantity $Q$ evolves according to
\begin{multline}\label{FirstMonotonicity}
\frac{d}{dt} Q(t)
=-\int_I \Lambda \left| K+\rho\frac{W^\perp}{2(t^*-t)} \right|^2_g ds 
+\left. \left\langle \left(\Lambda {\frac{\rho}{2(t^*-t)}} W \right) , T \right\rangle_g \right|_{\partial I} 
\\
- \int_I  \frac{\Lambda}{2(t^*-t)} \left(1- \left\langle T,W \right\rangle_g^2\right) \left(\rho\frac{\scrS'(\rho)}{\scrS(\rho)}-1\right) ds,
\end{multline}
where $W^\perp = \langle W,N\rangle_g N$.
\end{lemma}

The idea is to integrate \eqref{FirstMonotonicity}; first we consider the boundary terms and note that for $\tilde u\in \partial I$ we have
\begin{equation}\label{IntegralOfLambda}
\int_0^{t^*} \frac{1}{2(t^*-t)} \Lambda(\tilde u,t) dt <\infty,
\end{equation}
while $\left|\left\langle {\rho(\tilde u, t)} W  , T \right\rangle_g\right|$ is bounded uniformly in $t$.

We next consider the final integral in \eqref{FirstMonotonicity}, noting that $|\rho|$ and $1-\langle T,W\rangle^2_g $ are bounded.
As
\begin{equation}
\frac{ \scrS^\prime(\phi)}{\scrS} = \frac{1}{\phi} + \mathcal{O}(|\phi|^2)
\end{equation}
we see that 
$|\rho\frac{\scrS'(\rho)}{\scrS(\rho)} -1| = \mathcal{O}(\rho^3)$ is bounded as well.
Thus \eqref{IntegralOfLambda} implies that the time integral of the the final term of \eqref{FirstMonotonicity} is bounded.
Since $\lim_{t\to t^*}Q(t)$ is finite, we arrive at the following.
\begin{prop}
\begin{equation}\label{FirstHomotheticBound}
\int_0^{t^*}\int_I \Lambda \left| K+\rho\frac{W^\perp}{2(t^*-t)} \right|^2_g ds\,dt <\infty.
\end{equation}
\end{prop}

\section{Re-scalings of the curve shortening flow}\label{S-Rescalings}
Our strategy is to use the results of Sections \ref{S-DistanceComparison} and \ref{S-Monotonicity} to analyze potential singularities, in order to rule out the possibility that a singularity occurs.
In the analysis of singularities it is useful to consider re-scalings of the flow; compare the following to \cite{Angenent90I},\cite{AbreschLanger86},\cite{Altschuler91}.

Given a  function $R(t)>0$, we define a time-dependent re-scaled geometry on $\Sigma$ by $g_R = R^2g$.
Suppose $F:I\times[0,t^*)\to \Sigma$ solves \eqref{main} and let $T_R$, $N_R$, and $K_R$ be the unit tangent, unit normal, and curvature vectors determined by $g_R$.
We easily see that $T_R = \frac{1}{R} T$, $N_R= \frac{1}{R}N$, and $K_R = \frac{1}{R^2} K$.
Set $\kappa_R = \langle K_R, N_R\rangle_{g_R}$; note that $\kappa_R = \frac{1}{R}\kappa$.
We define a re-scaled time variable $\tau$ by $d\tau = R^2 dt$ and use the notation $\tau^0 = \tau(0)$, $\tau^* = \tau(t^*)$.
We see that $\tilde F:I\times[\tau_0,\tau^*)\to\Sigma$ given by $\tilde F(u, \tau) = F(u, t(\tau))$ satisfies $\frac{d}{d\tau} \tilde F = K_R$. 

Under rescalings of the geometry on $\Sigma$ by function $R$, we obtain re-scaled distance functions $D_R = RD$ and $L_R = RL$. 
Thus the ratio $D/L$ is scale-invariant and the bound in Proposition \ref{DistanceComparisonBound} holds for re-scaled flows as well.

\section{Singularities}\label{S-Singularities}
Proposition \ref{ContinuationCriterion} implies that either the flow exists for all time, or the curvature $\kappa$ is unbounded on some finite time interval.
We first address the latter case, assuming in this section that $[0,t^*)$, with $t^*<\infty$, is the maximal time interval on which a smooth flow exists and hence that $\sup_{[0,t^*)}\|\kappa\|_{L^\infty} = \infty$.

In light of \eqref{EvolveKappa}, the maximum principle implies a lower bound on the rate of curvature blowup. 
In fact the following estimate is true for a broad class of flows; see \cite{Angenent90I}.
\begin{lemma}\label{BlowupBoundedBelow}
The quantity $\sqrt{t^*-t}\,\|\kappa(t)\|_{L^\infty}$ is bounded below, strictly away from zero, for $t\in[0,t^*)$.
\end{lemma}
This leads to the following definition:
A singularity is called \emph{type 1} if 
$$\sup_{t\in[0,t^*)} \sqrt{t^*-t}\,\|\kappa(t)\|_{L^\infty} <\infty,$$ 
while the singularity is otherwise called \emph{type 2}.

The following is useful.
\begin{lemma}\label{ConstructBlowupSequence}
There exists sequence $(u_n,t_n)$ converging to $(u^*,t^*)$ such that
\begin{enumerate}
\item $|\kappa(u_n,t_n)|\to\infty$ as $n\to\infty$,
\item $\|\kappa(t)\|_{L^\infty}\leq |\kappa(u_n,t_n)|$ whenever $t\leq t_n$, and
\item there exists $P^*\in \Sigma$ with $F(u_m,t_n) \to P^*$.
\end{enumerate}
\end{lemma}
The construction is as in \cite{Altschuler91}.
Note that $\|\kappa(t)\|_{L^\infty}<\infty$ for $t<t^*$. 
Fix $\epsilon>0$ and set $E = \{ t\leq t^*\mid \|\kappa(t^\prime)\|_{L^\infty} \leq \|\kappa(t)\|_{L^\infty} \text{ for all } t^\prime\in[0,t]\}$; note $E$ is not empty.
Define the function $\sigma_\epsilon(t) = \min\{\sigma\in[t,t^*)\mid \|\kappa(\sigma)\|_{L^\infty} = (1+\epsilon)\|\kappa(t)\|_{L^\infty}\}$.
Since $t\in E$ implies $\sigma_\epsilon(t)\in E$, we may fix $t_0\in E$ and define $t_{n+1} = \sigma_\epsilon(t_n)$.
Let $u_n$ be such that $\kappa(u_n,t_n) = \|\kappa(t_n)\|_{L^\infty}$.
Since the points $F(u_n,t_n)$ are bounded in $\Sigma$ we may pass to a convergent subsequence to obtain the desired sequence, which we call a blowup sequence.

A number results have been obtained concerning singularities for curve shortening flow, see especially \cite{Gage90b},\cite{Grayson89},\cite{Angenent90I,Angenent91II}; for results concerning planar flows see \cite{GageHamilton86},\cite{Grayson87},\cite{Angenent91}.
These conclude that as $t\to t^*$ the flow converges towards a limiting curve which, except for a finite number of singularities, is regular.
Furthermore, as $t\to t^*$ the curvature diverges on a sequence of arcs which turn through angle $\pi$.
While stated for closed curves, the results are local in nature and extend to curves with fixed endpoints.

We record the following two propositions from \cite{Angenent90I,Angenent91II}, where larger classes of curvature flows are considered.
\begin{prop}
As $t\to t^*$, the sets $F(I,t)$ converge in the Hausdorff metric to a curve $\Gamma_*$, which is piecewise smooth and has only a finite number of singular points.
\end{prop}

In order to estimate the size of the blowup set, define for $\epsilon>0$ the quantity $\alpha_\epsilon(t)$ by
\begin{equation}\label{DefineAlphaEpsilon}
\alpha_\epsilon(t) = \sup_{|J|<\epsilon} \left|\int_{J} \kappa(u,t)ds\right|;
\end{equation}
here $J$ is a sub-interval of $I$ and $|J| = \int_Jds$. In \cite{Angenent91II}, this definition is extended to non-smooth curves.

\begin{prop}\label{AlphaBound}
If $t^*<\infty$ we have
$
\limsup_{t\to t^*}\alpha_\epsilon(t) \geq\pi
$
for all $\epsilon>0$.
\end{prop}
In fact Proposition \ref{AlphaBound} is true in a neighborhood of each singular point.
We outline a proof of this latter proposition as the construction is used again below.

Fixing $\epsilon>0$ and blowup sequence $(u_n,t_n)$, let $B\subset\Sigma$ be the geodesic ball having radius $\rho<\frac{\epsilon\theta}{2}$ and centered at $P^*$; without loss of generality we have $F(u_n,t_n)\in B$.
Here $\theta$ is as in Proposition \ref{DistanceComparisonBound}; thus for $v_1,v_2\in I$ with images in $B$ we have $\int_{[v_1,v_2]}ds <\epsilon$.

Define a sequence of re-scalings of the flow (see Section \ref{S-Rescalings}) by constant factors $R_n = |\kappa(u_n,t_n)|$. 
Thus for each $n$ we have a re-scaled function $\tilde{F}_n$, which satisfies $\frac{d}{d\tau} \tilde{F}_n = K_{R_n}$ and is defined for re-scaled time $\tau = R^2_n(t-t_n)$ in the interval $[\tau^0_n,\tau^*_n)= [-R_n^2t_n, R_n^2(t^*-t_n))$.

Restrict each $\tilde{F}_n$ to the largest intervals $J_n(\tau)\subset I$ containing $u_n$ such that $\tilde{F}_n(J_n(\tau),\tau)\subset B$ for $\tau\in[\tau^0_n,\tau^*_n)$.
Parametrize $J_n$ by arclength $s_n$, induced by $g_{R_n}$, with $s_n(u_n)=0$.
Clearly $J_n\to\bbR$ as $n\to\infty$ for all $\tau$.

Let $D_n\subset\bbR^2$ be the disk of radius $\rho R_n$.
Let $\psi_n:D_n\to B$ correspond to geodesic polar coordinates on with respect to the rescaled metric $g_{R_n}$ and set $g_n = \psi_n^*g_{R_n}$. 
Note that the $(D_n,g_n)$ converge to the flat Euclidean plane as $n\to\infty$.

Let $F_n = (\psi_n^{-1})\circ\tilde{F}_n$, $K_n =\psi_n^*K_{R_n}$, $N_n =\psi_n^*N_{R_n}$, $\kappa_n = \langle K_{n}, N_{n}\rangle_{g_n}$.
The function ${F}_n$ satisfies $\frac{d}{d\tau} F_n = K_n$ for $\tau \in [\tau^0_n,\tau^*_n)$.
Note $|\kappa_n(0,0)|=1$; thus by Lemma \ref{ConstructBlowupSequence} we have $|\kappa_n|\leq 1$ for $\tau\leq 0$; this implies bounds on higher derivatives via Lemma \ref{L-KappaBoundIsEnough} and leads to the following.
\begin{lemma}
Upon passing to a subsequence, the sequence $F_n$ converges to limiting flow $F_\infty:\bbR\times(-\infty,0]\to\bbR^2$ satisfying \eqref{main} and having curvature $\kappa_\infty$ bounded by $1$.
\end{lemma}
We remark that in the case of a type 2 singularity the flow $F_\infty$ is defined for all $\tau\in\bbR$ and has curvature bounded as $\tau\to\infty$; see Section \ref{S-TypeTwoSingularities} below.

We proceed to establish the proposition by contradiction, assuming that $\alpha_\epsilon(t)\leq\alpha_0<\pi$ for $t$ sufficiently close to $t^*$.
Since the integral in \eqref{DefineAlphaEpsilon} is scale-invariant, we have 
\begin{equation}
\left|\int_{[s_1,s_2]}\kappa_\infty(\tau)\,ds\right|\leq\alpha_0
\end{equation}
for any $s_1,s_2\in\bbR$ and $\tau\leq 0$, which in turn implies that, for short time intervals and in appropriately chosen Cartesian coordinates, $F_\infty$ can be written as a graph of a function $y(x,\tau)$, which must satisfy
\begin{equation}
\partial_\tau y 
= \frac{1}{(1+(\partial_xy)^2)^2}\partial^2_xy.
\end{equation}

Foliating the plane by straight lines, Proposition \ref{P:AngMax} implies that the gradient remains bounded and thus if $F_\infty(\cdot, \tau_1)$ is a graph, then the flow remains a graph in these coordinates for $\tau\in[\tau_1,0]$.
Taking $\tau_1\to-\infty$ we are able to represent $F_\infty$ as a graph for $\tau\in(-\infty,0]$.
In fact, as $\partial_xy$ is a bounded solution to a parabolic equation
the Harnack inequality of \cite{Moser64} implies that $\partial_xy$ is constant.
Since $|\kappa_\infty(0,0)|=1$ we obtain the desired contradiction and establish Proposition \ref{AlphaBound}.

\section{Type 1 singularities}\label{S-TypeOneSingularities}

Following Huisken \cite{Huisken90}, the monotonicity formula can be used to classify type 1 singularities as corresponding to self-similar solutions to planar curve shortening flow.
These self-similar solutions are either not embedded, do not have sufficient curvature, or do not satisfy the lower bound on $D/L$; thus we obtain the following.

\begin{prop}\label{Type1Flow}
Under our assumptions, type 1 singularities do not occur.
\end{prop}

Suppose $F:I\times [0,t^*)\to \Sigma$ has a type 1 singularity at $t^*$ and consider the re-scaled flow $\tilde{F}$ given by taking  $R(t)  = (2(t^*-t))^{-1/2}$.
The re-scaled flow is defined for $\tau\in [\tau_0,\infty)$.

Lemma \ref{BlowupBoundedBelow} and the definition of type 1 singularity implies that and the re-scaled curvature $\kappa_R$ is bounded above satisfies $c \leq |\kappa_R(u_n,t_n)|\leq C$ along any sequence satisfying (1) and (3) of Lemma \ref{ConstructBlowupSequence}.

Following the construction in Section \ref{S-Singularities}, we may restrict $F$ to a small ball $B$ about the singular point.
Under the re-scaling, $(B,g_R)$ converges to the flat Euclidean plane as $\tau \to \infty$ and, upon re-parametrization, the (restricted) re-scaled flow $\tilde F$ converges to a limiting curve $F_\infty:\bbR\to \bbR^2$.
The limiting curve is embedded, has infinite arc-length, is unbounded, and has bounded curvature; it cannot have vanishing curvature due to the lower near the singular point $P^*$.
In view of Proposition \ref{AlphaBound}, the curve turns through an angle of at least $\pi$.

Changing to re-scaled variables in \eqref{FirstHomotheticBound} yields
\begin{equation}\label{MainHomotheticBound}
\int_{\tau_0}^{\infty} R \int_I  \Lambda_R \left| K_R+\rho_RW^\perp_R \right|^2_{g_R} ds\,d\tau <\infty,
\end{equation}
where $\Lambda_R = \frac{1}{\sqrt{2\pi}}\exp{\left(-\frac12 \rho_R^2\right)}$.
Thus there exists a sequence of times $\tau_n \to \infty$ along which
\begin{equation}\label{WeakPointwiseConvergence}
\left.  \Lambda_R \left| K_R+\rho_RW^\perp_R \right|^2_{g_R} \right|_{\tau_n} \to 0
\end{equation}
for almost all points in $I$.
Therefore the limiting curve satisfies
\begin{equation}\label{Homothetic}
 K+F_\infty^\perp =0.
\end{equation}
In particular, the limiting curve corresponds to a solution to planar curve shortening flow moving by homothety.
Such curves are classified in \cite{AbreschLanger86},\cite{Halldorsson2010}.
The relevant solutions to \eqref{Homothetic} are either asymptotic to a cone and thus have total curvature strictly less than $\pi$, are not embedded, or do not satisfy any lower bound on $D/L$.
Thus we obtain a contradiction and conclude that a type 1 singularity cannot occur.

\section{Type 2 singularities}\label{S-TypeTwoSingularities}

We rule out the possibility that a type 2 singularity occurs by showing any such singularity is asymptotic, after re-scaling, to the planar `grim reaper' solution $y = t-\log{\cos{x}}$, which violates lower bound on $D/L$.

\begin{prop}\label{Type2Flow}
Under our assumptions, a singularity of type 2 does not occur.
\end{prop}

Our argument follows those in \cite{Hamilton95},\cite{Altschuler91}.
Let $F_n$ be the sequence of flows constructed in the proof of Proposition \ref{AlphaBound} and let $F_\infty:\bbR\times(-\infty,0]\to\bbR^2$ be the limit flow.

Since $\int|\kappa|ds$ is invariant under re-scaling, and by Lemma \ref{TurningAngleEstimate} is monotone decreasing, we see that $\int |\kappa_\infty|ds$ is finite and constant.
Thus, applying Lemma \ref{TurningAngleEstimate} to $F_\infty$ we have
\begin{equation}
\int_{-\infty}^0\sum_{\{u \mid \kappa(u)=0\}}|\partial_s\kappa(u)| d\tau  =0,
\end{equation} 
which implies inflection points of the limit curve are degenerate.
Proposition \ref{P:AngMax} implies that such a curve must be a straight line, which cannot be as $|\kappa_\infty(0,0)|=1$.
Thus $\kappa_\infty\neq0$; in view of Proposition \ref{AlphaBound}, we conclude that the limiting flow $F_\infty$ is convex at all $\tau$ and turns exactly through angle $\pi$.

We next claim that $F_\infty$ exists for all time with bounded curvature.

Let the blowup sequence constructed as in the proof of Lemma \ref{ConstructBlowupSequence}.
Note $\|\kappa(t_m)\|_{L^\infty} = (1+\epsilon)^{m-n} \|\kappa(t_n)\|_{L^\infty}$.
Thus
\begin{equation}
\begin{aligned}
(t^*-t_n)\|\kappa(t_n)\|_{L^\infty}^2 
&= \sum_{m=n}^\infty (t_{m+1}-t_m)\|\kappa(t_n)\|_{L^\infty}^2
\\
&\leq \left(\limsup_{m\to\infty}\|\kappa(t_m)\|_{L^\infty}^2(t_{m+1}-t_m) \right) \sum_{m=0}^\infty(1+\epsilon)^{-2m}
\end{aligned}
\end{equation}
and, since the sum is bounded, $\limsup_{m\to\infty}\|\kappa(t_m)\|_{L^\infty}^2(\sigma_\epsilon(t_{m})-t_m)=+\infty$.
Upon rescaling we have $\|\kappa_n(0)\|_{L^\infty} =1$.
Hence the interval $[0,\tau({\sigma}_\epsilon(t_n))]$, on which $\|\kappa_n\|_{L^\infty} \leq 1+\epsilon$, tends to $[0,\infty)$ as $n\to\infty$.
We conclude that the flow $F_\infty$ may be extended for $\tau\in (-\infty,\infty)$ with curvature bounded by $1+\epsilon$ for any $\epsilon>0$.

Not only is the curvature bounded, but $|\kappa_\infty|$ is also integrable. Thus the curvature tends to zero for large $s$.
Since the limiting flow $F_\infty$ is convex and eternal, results in \cite{Hamilton95b} state that $F_\infty$ must be a translating soliton.
By the classification in \cite{Halldorsson2010} such a curve is necessarily the grim reaper (see also \cite{Altschuler91}), for which no lower bound on $D/L$ is possible.
Thus we conclude that no singularity of type 2 forms.

\section{Convergence to geodesics}\label{S-Converge}
In light of Propositions  \ref{Type1Flow} and \ref{Type2Flow}, the continuation criterion in Proposition \ref{ContinuationCriterion} implies that the flow $F$ exists for $t\in[0,\infty)$.
We now show that the curvature tends to zero, which implies that there exist sequences of times $t_k\to\infty$ such that $F(I,t_k)$ tend to a geodesic.
When there exists a unique geodesic connecting $A$ and $B$, then this statement can be strengthened: The flow converges to that geodesic.
The lack of uniqueness in the case of conjugate points is somewhat expected; an analogous situation arises for closed curves, see \cite{Gage90b},\cite{Grayson89}.

\begin{prop}\label{P-CurvatureGoesToZero}
Suppose the a solution to \eqref{main}-\eqref{DirichletBC} exists for $t\in[0,\infty)$.
Then all derivatives of  $\kappa$ tend pointwise to zero.
\end{prop}

\begin{prop}\label{P-AchieveGeodesic}
If the endpoints $A,B$ lie in a strictly convex set (so that they are not conjugate), then the flow converges to the unique geodesic joining $A$ and $B$.
\end{prop}

We outline the argument establishing Proposition \ref{P-CurvatureGoesToZero}, following \cite{Gage90b},\cite{Grayson89}.

First we show that if the flow exists for all time, then $\int_I\kappa^2ds \to 0$ and $\int_I(\partial_s\kappa)^2\to 0$ as $t\to\infty$.
Note that by integrating from an endpoint, where $\kappa=0$, and using Cauchy-Schwartz we have
\begin{equation}\label{CheapSobolev}
\|\kappa\|_{L^\infty}^2 \leq \scrL\int_I(\partial_s\kappa)^2ds.
\end{equation}
Thus Lemma \ref{L-EvolutionOfCurvature} implies
\begin{equation}\label{EvolveL2}
\begin{aligned}
\frac{d}{dt}\int_I\kappa^2ds &= \int_I \left( -2(\partial_s\kappa)^2 + \kappa^4 + 2\kappa^2 S\right)ds
\\
&\leq \scrL\left(\int_I (\partial_s\kappa)^2ds\right)\left(\int_I \kappa^2ds - \frac{2}{\scrL}\right) + 2|S|\int_I\kappa^2ds.
\end{aligned}
\end{equation}
Note that $\scrL$ is bounded above by $\scrL(0)$ and below by the distance between the endpoints.

From Lemma \ref{L-evol} we have 
\begin{equation}\label{L1L2Curvature}
\int_0^\infty \int_I \kappa^2ds\, dt = \scrL(0) - \lim_{t\to\infty}\scrL(t)<\infty
\end{equation} 
and thus have $t_k\to\infty$ with 
\begin{equation}
\int_I\kappa(t_k)^2ds \to 0\quad\text{ and }\quad\int_{t_k}^\infty \int_I\kappa(t)^2ds\,dt\to 0.
\end{equation}

For small $\epsilon>0$, choose $t_k$ such that $\int_I\kappa(t_k)^2ds <  \frac12\epsilon$ and $\int_{t_k}^\infty \int_I\kappa^2ds\,dt<\epsilon^2$.
Suppose at $t^\prime>t_k$ we have $\int_I\kappa(t^\prime)^2ds = \epsilon$.

For $\epsilon$ small, depending on $\scrL(0)$ the first term in \eqref{EvolveL2} is negative and can be dropped.
Integrating from $t_k$ to $t^\prime$ we find
\begin{equation}
\frac12\epsilon<\int_I\kappa(t^\prime)^2ds - \int_I\kappa(t_k)^2ds 
\leq 2|S|\int_{t_k}^{t^\prime}\int_I\kappa^2ds\,dt
<2|S|\epsilon^2,
\end{equation}
which is a contradiction for small $\epsilon$.
Thus for any $\epsilon>0$ there exists $t_k$ such that $\|\kappa(t)\|_{L^2}^2\leq \epsilon$ for $t>t_k$ and the convergence of $\|\kappa\|_{L^2}$ is proved.

A similar argument shows the derivative of curvature tends to zero in $L^2$.
Integrating from the point of maximum curvature we have
\begin{equation}\label{SecondCheapSobolev}
\|\partial_s\kappa\|_{L^\infty}^2 \leq \scrL\int_I(\partial_s^2\kappa)^2ds.
\end{equation}
Thus we compute
\begin{equation}\label{DerivativeOfL2Kappa}
\begin{aligned}
\frac{d}{dt}\int_I(\partial_s\kappa)^2ds
&= \int_I \left( 7\kappa^2(\partial_s\kappa)^2-2(\partial^2_s\kappa)^2 - 2S\kappa (\partial^2_s\kappa)  \right)ds
\\
&\leq \left(7\scrL \int_I\kappa^2ds -1\right) \int_I(\partial^2_s\kappa)^2ds
	+ |S|^2 \int_I\kappa^2ds.
\end{aligned}
\end{equation}
The first term is negative for large $t$; thus \eqref{L1L2Curvature} ensures that $\int_I(\partial_s\kappa)^2ds\to 0$.

Similar arguments using \eqref{HowDerivativesEvolve} can be used to conclude that $\|\partial_s^k\kappa\|_{L^2}\to 0$ for all $k$; see \cite{Gage90b},\cite{Grayson89} for details.
The Sobolev theorem yields pointwise estimates and we conclude that $\|\partial_s^k\kappa\|_{L^\infty}\to 0$ for all $k$.

Proposition \ref{P-CurvatureGoesToZero} is sufficient to conclude, via the Arzela-Ascoli theorem, that some subsequence of curves $F(I,t_k)$ converges to a geodesic.
In order to establish the uniqueness of the limiting curve, and hence Proposition \ref{P-AchieveGeodesic}, we show that $\|\kappa(t)\|_{L^\infty}$ is integrable in time.

We make use of Wirtinger's inequality (see \cite{GageHamilton86}):
If $g(a) = 0=g(b)$ with $b-a\leq\pi$, then
\begin{equation}
\int_a^b g(\eta)^2d\eta \leq \int_a^b (g^\prime(\eta))^2d\eta.
\end{equation}
An immediate consequence, with $d\eta = (\pi/\scrL)ds$, is
\begin{equation}\label{Wirtinger}
\int_I\kappa^2ds \leq \left(\frac{\scrL}{\pi}\right)^2\int_I (\partial_s\kappa)^2ds.
\end{equation} 

In the case of curves on $\bbS^2$, the assumption that endpoints $A,B$ lie in a strictly convex set, implies that along the convergent sequence we have $\scrL(t_k)$ tending to some limit strictly less than $\pi$.
Since $\scrL$ is monotone decreasing, we may assume that $\scrL<\pi$ whenever $S=1$.
Thus, computing as in \eqref{EvolveL2}, we have
\begin{equation}
\begin{aligned}
\frac{d}{dt}\int_I e^{\delta t}\kappa^2ds
&= e^{\delta t}\int_I\left(\delta\kappa^2 -2(\partial_s\kappa)^2 + \kappa^4 + 2S\kappa^2 \right)ds
\\
&\leq e^{\delta t} \left(\delta\frac{\scrL^2}{\pi^2} + 2\frac{\scrL^2}{\pi^2}\max(0,S) + \|\kappa\|^2_{L^\infty} -2\right)\int_I(\partial_s\kappa)^2ds.
\end{aligned}
\end{equation}
Since either $S\leq 0$ or $\scrL<\pi$ we may find $\delta>0$ such that the right side is negative and
\begin{equation}
\int_I \kappa^2ds \leq Ce^{-\delta t}
\end{equation}
for some constant $C\gg 1$.

In order to control the first derivative, consider those times when $C\int_I\kappa^2ds <\int_I(\partial_s\kappa)^2ds$.
For these times we have, by the H\"older inequality, that
\begin{equation}
C\int_I\kappa^2ds <\int_I(\partial_s\kappa)^2ds 
= \int_I \kappa (\partial_s^2\kappa)ds 
\leq \left(\int_I\kappa^2ds\right)^{1/2}\left( \int_I(\partial_s^2\kappa)^2\right)^{1/2}.
\end{equation}
Thus \eqref{DerivativeOfL2Kappa} yields
\begin{equation}
\frac{d}{dt}\int_I(\partial_s\kappa)^2ds \leq \left(7\scrL Ce^{-\delta t} -1 + \frac{|S|}{C^2} \right)\int_I(\partial_s^2\kappa)^2ds.
\end{equation}
Since \eqref{SecondCheapSobolev} implies $\int_I (\partial_s\kappa)^2ds \leq \scrL^2 \int_I(\partial_s^2\kappa)^2ds$ we may choose $t$ large enough so that
\begin{equation}
\frac{d}{dt}\int_I(\partial_s\kappa)^2ds \leq - \frac{1}{2\scrL(\infty)^2}\int_I(\partial_s\kappa)^2ds.
\end{equation}
Thus whenever $C\int_I\kappa^2ds <\int_I(\partial_s\kappa)^2ds$ we see that $\int_I(\partial_s\kappa)^2ds$ is exponentially decaying.

From \eqref{CheapSobolev} we see that $\|\kappa(t)\|_{L^\infty}$ decays exponentially, and thus is integrable.
By integrating the flow, we obtain a unique limit curve which, by Proposition \ref{P-CurvatureGoesToZero} is smooth and has zero curvature.

\bibliographystyle{plain}	
\bibliography{ALT}
\end{document}